\documentclass[2p]{article}
\usepackage{amssymb,amsmath,amsthm}
\usepackage{indentfirst}
\usepackage{exscale}
\usepackage{relsize}
\usepackage[numbers,sort&compress]{natbib}

\usepackage{geometry}
\usepackage{color}
\newcommand{\R}{\mathbb{R}}
\theoremstyle{definition}

\allowdisplaybreaks
\theoremstyle{remark}

\numberwithin{equation}{section}

 \UseRawInputEncoding

\begin{document}

\title{\Large\bf{ Infinitely many solutions for a biharmonic-Kirchhoff system on locally finite graphs }
 }
\date{}

\author {Xiaoyu Wang$^{1}$, Junping Xie$^{2}$ \footnote{Corresponding author, E-mail address: hnxiejunping@163.com}, \ Xingyong Zhang$^{1,3}$\\
{\footnotesize $^1$Faculty of Science, Kunming University of Science and Technology, Kunming, Yunnan, 650500, P.R. China.}\\
{\footnotesize $^2$Faculty of Transportation Engineering, Kunming University of Science and Technology,}\\
 {\footnotesize Kunming, Yunnan, 650500, P.R. China.}\\
{\footnotesize $^{3}$Research Center for Mathematics and Interdisciplinary Sciences, Kunming University of Science and Technology,}\\
 {\footnotesize Kunming, Yunnan, 650500, P.R. China.}\\}

 \date{}
 \maketitle

 \begin{center}
 \begin{minipage}{15cm}

\small  {\bf Abstract:} The study on the partial differential equations (systems) in the graph setting is a hot topic in recent years because of their applications to image processing and data clustering. Our motivation is to develop some existence results for biharmonic-Kirchhoff systems and biharmonic systems in the Euclidean setting, which are the continuous models, to the corresponding systems in the locally finite graph setting, which are the discrete models. We mainly focus on the existence of infinitely many solutions for a biharmonic-Kirchhoff system on a locally finite graph. The method is variational and the main tool is the symmetric mountain pass theorem.  We obtain that the system has infinitely many solutions when the nonlinear term admits the super-$4$ linear growth, and we also present the corresponding results to the biharmonic system. We also find that the results in the locally finite graph setting are better than that in the Euclidean setting, which caused by the better embedding theorem in the locally finite graph.

\end{minipage}
 \end{center}
  \allowdisplaybreaks
 \vskip2mm
\section{Introduction and main results}\label{section 1}
In this paper, we investigate the existence of infinitely many solutions for the following biharmonic-Kirchhoff system on a locally finite graph $G=(V,E)$:
\begin{eqnarray}
\label{101}
   \begin{cases}
 \Delta^2 u- (a_1+b_1 \int_V |\nabla u|^2 d\mu)\Delta u + V_1(x)u=F_u(x,u,v),\;\;\;\;\hfill x\in V,\\
 \Delta^2 v- (a_2+b_2 \int_V |\nabla v|^2 d\mu)\Delta v + V_2(x)v=F_v(x,u,v),\;\;\;\;\hfill x\in V,\\
   \end{cases}
\end{eqnarray}
where $V$ is the vertexes set of $G$, $E$ is the edges set of $G$, $\Delta^2:=\Delta(\Delta)$ is the biharmonic operator, $a_i>0,\;b_i\ge 0$, $V_i:V\to \R$, $i=1,2$, $F:V\times \R\times\R \to \R$, $F_u(x,u,v)=\frac{\partial F(x,u,v)}{\partial u}$ and $F_v(x,u,v)=\frac{\partial F(x,u,v)}{\partial v}$.
\par
Next, we provide the precise explanations on the notations in system (\ref{101}). We refer to \cite{grigoryan2016, han2020} for these notations.  $xy\in E$ denotes the edge which links $x$ to its neighbor $y$ in $V$ and denote $\omega_{xy}>0$ by its weight. For any edge $xy\in E$, it is supposed that $\omega_{xy}=\omega_{yx}$. The graph $G=(V,E)$ is also supposed to be a connected graph which means that any two vertices $x$ and $y$ in $V$ have to be linked by finite edges.
Moreover, $G=(V,E)$ is a locally finite graph which means that the number of edges $xy \in E$ is finite for any given $x \in V$.
 $dist(x,y)$ denotes the distance of two vertices $x,y$, which means the minimal number of edges that link $x$ with $y$.
 Let $C(V)=\{u|u:V\to\R\}$ and define the Laplacian operator $\Delta: C(V)\to C(V)$ by
\begin{eqnarray*}
\label{no001}
\Delta u(x)=\frac{1}{\mu(x)}\sum_{y\thicksim x}\omega_{xy}(u(y)-u(x)),
\end{eqnarray*}
where $\mu:V\to \R^+$ is a finite measure with positive lower bound $\mu_{\min}$, that is, $\mu(x)\geq \mu_{\min}>0$, and $y\thicksim x$ denotes the vertex $y$ satisfying the edge $xy\in E$. Then
\begin{eqnarray*}
\label{no001}
\Delta^2 u(x):=\Delta (\Delta u(x))=\frac{1}{\mu(x)}\sum_{y\thicksim x}\omega_{xy}(\Delta u(y)-\Delta u(x)).
\end{eqnarray*}
Define the gradient $\nabla u$   of $u$ at $x \in V$ by a vector
\begin{eqnarray*}
\nabla u (x):=\left(\sqrt{\frac{w_{xy}}{{2\mu(x)}}}(u(y)-u(x))\right)_{y\thicksim x, y\in V}
\end{eqnarray*}
which is indexed by $y$. Let
\begin{eqnarray*}
\label{no002}
\Gamma(u,v)(x):=\nabla u(x)\cdot \nabla v(x)=\frac{1}{2\mu(x)}\sum\limits_{y\thicksim x}w_{xy}(u(y)-u(x))(v(y)-v(x)).
\end{eqnarray*}
Then
\begin{eqnarray*}
\label{no003}
|\nabla u|(x)=\sqrt{\Gamma(u,u)(x)}=\left(\frac{1}{2\mu(x)}\sum\limits_{y\thicksim x}w_{xy}(u(y)-u(x))^2\right)^{\frac{1}{2}}.
\end{eqnarray*}
For any $u\in C(V)$, we denote
\begin{eqnarray*}
\label{no004}
\int\limits_V u(x) d\mu=\sum\limits_{x\in V}u(x)\mu(x).
\end{eqnarray*}
The biharmonic operator $\Delta ^2 u$ of $u : V \rightarrow \R$ satisfies the following identity (see \cite{han2020}):
\begin{eqnarray*}
\label{no005}
\int_V (\Delta^2 u)\phi d\mu =\int_V \Delta u \Delta \phi d\mu,\; \forall \phi \in  C_c(V),
\end{eqnarray*}
where $C_c(V):= \{u \in C(V)| \mbox{the set }\{x \in V : u(x) \neq 0\} \;\;\text{is of finite cardinality}\}$.
\par
In recent years, the analysis on locally finite graph has been becoming a hot topic because of its own mathematical interest and its applications to image processing and data clustering (see \cite{elmoataz2012, elmoataz2015, elmoataz20171, elmoataz20172, ennaji2023}). Especially, the development of the analysis on locally finite graph greatly promote the investigation of the partial differential equations on locally finite graph (see \cite{grigoryan2016, han2020, han2021, pan2023, pang2023, shao2023, shao2024, wang2024, yang2024, zhang2024, yu2024, ou2023}), which will probably benefit to the numerical computation of the partial differential equations. It is well known that variational methods and critical point theory are the important tools to deal with
the existence and multiplicity of nontrivial solutions of differential equations or difference equations (for example, see \cite{gao2023, liu2012, shao2024, wei2012, yang2024, yu2024, yuan2024, zhao2016}). Naturally, it will come to mind whether one can apply these two tools to the partial differential equations on locally finite graphs. Indeed, there have been some interest works in this field. For example, in \cite{grigoryan2016},  Grigor¡¯yan-Lin-Yang used the mountain pass theorem  to prove that there is a positive solution to a  Laplacian equation with a Dirichlet boundary. They also applied the method to $p$-Laplacian equations and $poly$-Laplace equations. In \cite{han2021}, Han-Shao used the mountain pass theorem and Nehari manifold, under suitable conditions, it is proved that when $\lambda > 1$, there is a positive solution to a $p$-Laplacian equation with a parameter $\lambda$. And they proved the convergence of solutions when $\lambda \rightarrow \infty$. In \cite{pang2023}, Pang-Zhang studied the existence of nontrivial solutions of poly-Laplace systems with parameters on finite graphs and $(p,q)$-Laplace systems with parameters on locally finite graphs. The three solutions of systems were obtained by using a critical point theorem without compactness condition. In \cite{yang2024}, Yang-Zhang generalized two embedding theorems and investigated the existence and multiplicity of nontrivial solutions for a $(p,q)$-Laplacian coupled system with perturbations and two parameters $\lambda_1$ and $\lambda_2$ on locally finite graph. By using mountain pass theorem and Ekeland¡¯s variational principle, they obtained that system has at least one nontrivial solution when the nonlinear term satisfies the sub-$(p,q)$ conditions and obtain that system has at least one solution of positive energy and one solution of negative energy when the nonlinear term satisfies the super-$(p,q)$ conditions which is weaker than the well known Ambrosetti-Rabinowitz condition. In \cite{pan2023}, Pan-Ji use the constrained variational method to prove the existence of a least energy sign-changing solution $u_b$ of a Kirchhoff equation if the nonlinear term $f$ satisfy certain assumptions, and to show the energy of $u_b$ is strictly larger than twice that of the least energy solutions. In \cite{yu2024}, by using the well-known mountain pass theorem and Ekeland¡¯s variational principle, Yu-Xie-Zhang prove that there exist at least two fully nontrivial solutions for a $(p,q)$-Kirchhoff elliptic system with convex-concave nonlinearity and  perturbation terms on a locally finite graph. They also present a necessary condition of the existence of semitrivial solutions for the system. In \cite{ou2023}, Ou-Zhang use a direct non-Nehari manifold method to obtain the existence results of least energy sign-changing solutions and ground state solutions for a class of Kirchhoff-type equations with general power law, logarithmic nonlinearity and Dirichlet boundary value on a locally finite graph, and they obtain the sign-changing least energy is strictly larger than twice of the ground state energy. More results can be referred to \cite{bu2023, li2020, shuai2015, wen2019, ye2012, zhang2011, zhang2014, han2020}.
\par
In this paper, our works are mainly inspired by a recent work due to Han-Shao-Zhao \cite{han2020} where the following nonlinear biharmonic equation was investigated:
\begin{eqnarray*}
\label{106}
 \Delta^2 u- \Delta u + (\lambda a +1)u=|u|^{p-2}u,
\end{eqnarray*}
on a locally graph $G=(V,E)$. Under some suitable assumptions of Sobolev spaces defined, they proved that the equation admits a ground state solution $u_\lambda$ for any $\lambda>1$ and $p>2$. Moreover, they also demonstrated an asymptotic result of solutions $\{u_{\lambda}\}$  as $\lambda\to +\infty$. We are also inspired by the works due to Zhang-Tang-Zhang \cite{zhang2014} and Yuan-Liu \cite{yuan2024}. In  \cite{zhang2014}, Zhang-Tang-Zhang studied the following fourth-order elliptic equation in the Euclidean setting:
\begin{eqnarray}
\label{104}
   \begin{cases}
 \Delta^2 u- \Delta u + V(x)u=f(x,u),\;\;\;\;\hfill x\in \R^N,\\
 u \in H^2(\R^N),\\
   \end{cases}
\end{eqnarray}
where the potential $V \in C(\R^N,\R)$ is allowed to be sign-changing. Under the weak superquadratic conditions, they establish the existence of infinitely many solutions via variational method and symmetry mountain pass theorem for (\ref{104}).
In \cite{yuan2024} Yuan-Liu considered the following biharmonic-Kirchhoff equation in the Euclidean setting:
\begin{eqnarray}
\label{105}
   \begin{cases}
 \Delta^2 u- \left(a+b\int_{\R^5} |\nabla u|^2 dx \right)\Delta u + V(x)u=f(u),\\
 u(x)=u(|x|) \in H^2(\R^5),\\
   \end{cases}
\end{eqnarray}
where $V \in C(\R^5,\R)$ and $f \in C(\R,\R)$. By using a perturbation approach and the symmetric mountain pass theorem, the existence and multiplicity of radially symmetric  solutions for (\ref{105}) are obtained.

\par
Based on the works in \cite{han2020}, \cite{zhang2014} and \cite{yuan2024}, it is a natural motivation to apply the  variational method and symmetry mountain pass theorem to investigate the existence of infinitely many solutions for system (\ref{101}) which can be seen as a discrete analogy on locally finite graph for equations (\ref{104}) and (\ref{105}). As we will see in Remark 4.1 in section 4, there are some differences in the restrictions to the nonlinear term which are essentially caused by the better embedding relations between working spaces on locally finite graph than that on Euclidean space $\R^N(N\ge 2\mbox{ is an interger})$. Moreover, compared to equations (\ref{104}) and (\ref{105}) which are the single equation case, system (\ref{101}) possesses the difficulties caused by the coupling by $u$ and $v$ in nonlinear term so that we have to make more delicate analysis. Finally, in section 4, we also present some similar results for the  biharmonic system $(b_i=0,i=1,2)$ on locally finite graph corresponding to system (\ref{101}).
We firstly make the following assumptions:
\vskip2mm
\noindent
$(V)$ $V_i \in C(V, \R)$, $i=1,2$, $V_0:=\inf_V V_i(x)>0$ for all $x \in V$, there exists a vertex $x_0 \in V$ such that $V_i(x)\rightarrow +\infty$ as $d(x,x_0)\rightarrow +\infty$;\\
$(F_0)$ $F(x,0,0)\equiv 0$ for all $x\in V$, and $F(x,s,t)$ is continuously differentiable in $(t,s) \in \R^2$ for all $x \in V$;\\
$(F_1)$ there exists a function $a\in C(\R^+,\R^+)$ and a function $b:V\rightarrow \R^+$ with $b \in L^1(V)$ such that
$$|F_s(x,s,t)|,|F_t(x,s,t)|,|F(x,s,t)| \leq a(|(s,t)|)b(x)$$
for all $x\in V$,\;$(s,t)\in \R^2$;\\
$(C_1)$  $\lim\limits_{|(s,t)|\rightarrow 0} \frac{a(|(s,t)|)}{|(s,t)|^2}=h<+\infty$;\\
$(C_2)$ there exists a $\rho_*>0$ such that $\alpha_*:=\frac{1}{4}\rho_*^2 -\max\limits_{|s|+|t|\leq(\gamma_{\infty,1}+\gamma_{\infty,2})\rho_*} a(|(s,t)|) \int_V b(x)d\mu >0 $, where $\gamma_{\infty,i}=\left(\frac{1}{\mu_{\min} \inf\limits_V V_i(x)}\right)^{\frac{1}{2}}, i=1,2$;\\
$(F_2)$ $\lim_{|(s,t)|\rightarrow\infty} \frac{F(x,s,t)}{|(s,t)|^4}=+\infty$ uniformly for all $x \in V$;\\
$(F_3)$ there exists $r_0> 0$ such that
$$
F(x,s,t)\geq 0\;\text{for all}\;x \in V,\;(s,t) \in \R^2,\;|(s,t)| \geq r_0;
$$
$(F_4)$ $\mathcal{F}(x,s,t):= \frac{1}{4}\left(F_s(x,s,t)s+F_t(x,s,t)t\right)-F(x,s,t)\geq 0$, and there exist $c>0$ and $k > 1$ such that
$$|F(x,s,t)|^k \leq c |(s,t)|^{2k}\mathcal{F}(x,s,t)\;\text{for all}\;x \in V,\;(s,t) \in \R^2,\;|(s,t)| \geq r_0;$$
$(F_5)$ there exist $\mu > 4$ and $\sigma\in \left(0, \frac{\mu-2}{2\max\{\gamma_{2,1}^2,\gamma_{2,2}^2\}}\right)$ such that
$$\mu F(x,s,t) \leq sF_s(x,s,t)+tF_t(x,s,t)+\sigma (s^2+t^2)\;\text{for all}\;x \in V,\;(s,t) \in \R^2,$$
where $\gamma_{2,i} = \mu_{\min}(\inf\limits_V V_i(x))^{\frac{1}{2}} $;\\
$(F_6)$ $F(x,s,t)=F(x,-s,-t)$ for all $x \in V$ and $(s,t) \in \R^2$.

\vskip2mm
\par
Next, we state our main results in this paper.
\vskip2mm
\noindent
{\bf Theorem 1.1.} {\it Suppose that $b_i>0, i=1,2$, $(V)$, $(F_0), (F_1), (C_1) (\mbox{or } (C_2)), (F_2), (F_3), (F_4)$ and $(F_6)$ are satisfied. Then problem (\ref{101}) has infinitely many solutions $\{(u_k,v_k)\}_{k=1}^\infty$ such that
\begin{eqnarray}
\label{102}
&&\frac{1}{2} \int_V (|\Delta u_k|^2+ a_1|\nabla u_k|^2+ V_1(x)u_k^2)d\mu +\frac{1}{2} \int_V (|\Delta v_k|^2+ a_2|\nabla v_k|^2+ V_2(x)v_k^2)d\mu\nonumber\\
&&+\frac{b_1}{4}\left(\int_V |\nabla u_k|^2d\mu\right)^2+\frac{b_2}{4}\left(\int_V |\nabla v_k|^2d\mu\right)^2 - \int_V F(x,u_k,v_k)d\mu \rightarrow + \infty, \;\;\text{as} \;k \rightarrow + \infty .
 \end{eqnarray}}

\noindent
{\bf Theorem 1.2.} {\it Suppose that $b_i>0, i=1,2$, $(V)$, $(F_0), (F_1), (C_1) (\mbox{or } (C_2)), (F_2), (F_3), (F_5)$ and $(F_6)$ are satisfied.  Then problem (\ref{101}) has infinitely many solutions $\{(u_k,v_k)\}_{k=1}^\infty$  such that (\ref{102}) holds.}

\vskip2mm
\noindent
{\bf Remark 1.1} {\it } There are examples satisfying our assumptions in Theorem 1.1 and Theorem 1.2. For example, let
$$F_u(x,u,v)=a_1(x)(5u^5-3u^3),$$
$$F_v(x,u,v)=a_2(x)(5v^5-3v^3),$$
where $a_1(x),a_2(x)\in C(V,\R)$, $0< \inf_V a_1(x) \leq \sup_V a_1(x)< 1 $ and $0< \inf_V a_2(x)\leq \sup_V a_2(x)< 1$.
\vskip2mm
\noindent
{\bf Remark 1.2.} There are few papers to investigate the existence of infinitely many solutions for the partial differential equations on the locally finite graph except for \cite{pang2023} and \cite{yu2024} where the poly-Laplacian systems with parameters and Dirichlet boudary value and $(p,q)$-Laplacian systems were investigated and the methods and the restrictions for the nonlinear term $F$ are different from our Theorem 1.1 and Theorem 1.2. For the details, we refer readers to \cite{pang2023} and \cite{yu2024}.

 \section{Variational setting}\label{section 2}

\par
Define $L^r(V)=\{u:V \to \R \big|\int_V |u|^r d\mu <+\infty\}$ ($1\leq r<+\infty$) with the norm defined by
$$
\|u\|_r=
\left(\int_V|u(x)|^r d \mu\right)^{\frac{1}{r}}.
$$
Then  $(L^r(V),\|\cdot\|_r)$ is a Banach space  ($1\le r<+\infty$).
Define $L^\infty(V)=\left\{u:V \to \R \Big| \sup_{x\in V}|u(x)|<+\infty\right\}$  with the norm defined by
$$
\|u\|_{\infty}=\sup_{x\in V}|u(x)|.
$$
Let $W^{2,2}(V)$ be the completion of $C_c(V)$ under the norm
$$\|u\|_{W^{2,2}(V)}=\left( \int_V (|\Delta u|^2 + |\nabla u|^2 +u^2)d\mu \right)^{\frac{1}{2}}.$$
$W^{2,2}(V)$ is a Hilbert space with the inner product
$$(u,v)_{W^{2,2}(V)}= \int_V (\Delta u \Delta v + \nabla u \nabla v + uv)d\mu, \; u,v \in W^{2,2}(V).$$
We work in the following subspace of $W^{2,2}(V)$:
$$E_i :=\{u \in W^{2,2}(V): \int_V (|\Delta u|^2 + a_i|\nabla u|^2 +V_i(x)u^2)d\mu < +\infty\},\;i=1,2, $$
equipped with the inner product
$$(u,v)_{E_i}= \int_V (\Delta u \Delta v +a_i \nabla u \nabla v +V_i(x)uv)d\mu, \;u,v \in E_i,\;i=1,2 $$
and the norm
$$\|u\|_{E_i}=\left( \int_V (|\Delta u|^2 + a_i|\nabla u|^2 +V_i(x)u^2)d\mu \right)^{\frac{1}{2}},\; u,v \in E_i,\;i=1,2.$$
\vskip2mm
\noindent
{\bf Lemma 2.1.} {\it Assume that $(V)$ holds. Then, $E_i, i=1,2$ are continuously embedded into $L^r(V)$ for $ 2\le r\le+\infty$ and
\begin{eqnarray}
\label{221}
&  & \|u\|_r \leq\gamma_{r,i} \|u\|_{E_i},\;\forall u \in E_i,\;i=1,2,\ \  2\le r<+\infty,\\
\label{221a}&  & \|u\|_\infty\leq \gamma_{\infty,i}\|u\|_{E_i},\;\forall u \in E_i,\;i=1,2,
\end{eqnarray}
where $\gamma_{r,i} = \left(\mu_{\min}\inf\limits_V V_i(x) \right)^{\frac{2-r}{2r}}(\inf\limits_V V_i(x))^{-\frac{1}{r}} $ if $2\le r<+\infty$ and $\gamma_{\infty,i} = \left(\frac{1}{\mu_{\min} \inf\limits_V V_i(x)}\right)^{\frac{1}{2}}$ if $r=\infty$,  and the embedding $E_i \hookrightarrow L^s(V), i=1,2$ are compact for any $ 2 \leq r \leq +\infty$, that is, for any bounded sequence $\{u_k\}\subset E_i$, there exists $u \in E_i,i=1,2$ such that, up to a subsequence,
\begin{eqnarray*}
 \begin{cases}
\label{no007}
 u_k \rightharpoonup u \;\; \text{in}\; E_i,i=1,2;\\
  u_k(x)\rightarrow u(x) \;\; \forall x \in V;\\
  u_k\rightarrow u \;\;\  \text{in} \;L^r(V).
\end{cases}
\end{eqnarray*}}
\noindent
{\bf Proof.} This proof is similar to Lemma 2.5 in \cite{han2020} with some differences caused by $V_i(x),\;i=1,2$. For any given $ u \in E_i,\;i=1,2$ and any given vertex $x_0 \in V$, we have
\begin{eqnarray*}
\label{no008}
 \|u\|^2_{E_i} &=& \int_V (|\Delta u|^2 + a_i|\nabla u|^2 +V_i(x)u^2)d\mu   \\
\notag
& \geq & \int_V V_i(x)u^2 d\mu \\
\notag
&= &  \sum\limits_{x \in V} \mu(x) V_i(x)u^2(x) \\
\notag
& \geq &  \mu_{\min}  V_i(x_0)u^2(x_0)\\
\notag
& \geq &  \mu_{\min}  \inf\limits_V V_i(x)u^2(x_0),
\end{eqnarray*}
which implies
$$u(x_0)\leq \left(\frac{1}{\mu_{\min} \inf\limits_V V_i(x)}\right)^{\frac{1}{2}}\|u\|_{E_i}.$$
It follows from the arbitrariness of $x_0 \in V$ that $E_i \hookrightarrow L^\infty(V)$ is continuously embedding and (\ref{221a}) holds. Furthermore, we can get $E_i \hookrightarrow L^r(V)$ continuously for any $2\leq r< \infty$. In fact, for any $u \in E_i$, by (V), it is easy to see that $u \in L^2(V)$. Then, for any $2\leq s <\infty$,
$$
\int_V |u|^r d\mu= \int_V |u|^2|u|^{r-2} d\mu \leq (\mu_{\min}\inf\limits_V V_i(x) )^{\frac{2-r}{2}}\|u\|^{r-2}_{E_i} \int_V |u|^2 d\mu < +\infty,
$$
which implies that $u \in L^r(V)$ for any $2\leq r<\infty$ and
$$\|u\|_r =\left( \int_V |u|^r d\mu \right)^{\frac{1}{r}} \leq \left(\mu_{\min}\inf\limits_V V_i(x) \right)^{\frac{2-r}{2r}}(\inf\limits_V V_i(x))^{-\frac{1}{r}} \|u\|_{E_i} =\gamma_{r,i}\|u\|_{E_i}.$$
\par
Since $E_i(i=1,2)$ are Hilbert spaces, they are reflexive. Thus for any bounded sequence $\{u_{k}\}$ in $E_i$, there is, up to a subsequence, $u_{k}\rightharpoonup u_i$ in $E_i,\;i=1,2$. On the other hand, $\{u_{k}\}\subset E_i$ is also bounded in $L^2(V)$ and we also get $u_{k}\rightharpoonup u_i$ in $L^2(V)$, which tells us that, for any $\phi \in L^2(V)$,
\begin{eqnarray}
\label{201}
\lim\limits_{k\rightarrow \infty} \int_V (u_{k}-u_i)\phi d\mu=\lim\limits_{k\rightarrow \infty} \sum\limits_{x \in V} \mu(x)(u_{k}(x)-u_i(x))\phi(x)=0.
\end{eqnarray}
For any fixed $x_0 \in V$, let
\begin{eqnarray*}
\phi_0(x)=
 \begin{cases}
\label{no009}
1 \;\;\; x=x_0,\\
0 \;\;\; x\neq x_0.
\end{cases}
\end{eqnarray*}
Obviously, $\phi_0$ belongs to $L^2(V)$. By substituting $\phi_0$ into $\phi$ in (\ref{201}), we get
$$\lim\limits_{k\rightarrow \infty} \mu(x_0)(u_{k}(x_0)-u(x_0))=0,$$
which implies that
\begin{eqnarray}
\label{220}
\lim\limits_{k\rightarrow \infty} u_{k}(x)=u(x)\;\;\text{for any}\;\;x \in V.
\end{eqnarray}
\par
Next we prove $u_k \rightarrow u $ in $L^s(V)$ for all $2\leq s \leq \infty$. Since $\{u_{k}\}$ are bounded in $E_i,\;i=1,2$ and $u \in E_i,\;i=1,2$, there exist positive constants $C_i,\;i=1,2$ such that
$\|u_{k}-u\|^2_{E_i} \leq C_i,\;i=1,2$. Let $x_0 \in V$ be fixed. By $(V)$, for any $\epsilon > 0$, there exists some $r_0 >0$ such that
$V_i(x)\geq \frac{C_i}{\epsilon}$ as $d(x,x_0)> r_0$.
Then
\begin{eqnarray}
\label{203}
\int\limits_{d(x,x_0)>r_0} |u_k-u|^2 d\mu &\leq & \frac{\epsilon}{C_i} \int\limits_{d(x,x_0)>r_0} V_i(x)|u_k-u|^2 d\mu \\
\notag
&\leq & \frac{\epsilon}{C_i}\|u_k-u\|^2_{E_i}\\
\notag
&\leq & \epsilon.
\end{eqnarray}
Moreover, by (\ref{220}), we have
\begin{eqnarray}
\label{204}
\lim\limits_{k\rightarrow \infty} \int\limits_{d(x,x_0)\leq r_0} |u_k-u|^2 d\mu=0.
\end{eqnarray}
Combining (\ref{203}) and (\ref{204}), we conclude that
\begin{eqnarray*}
\label{no010}
\lim\limits_{k\rightarrow \infty} \int_V |u_k-u|^2 d\mu=0,
\end{eqnarray*}
that is, $u_k\rightarrow u$ in $L^2(V)$. Since
\begin{eqnarray*}
\label{no012}
\|u_k-u\|^2_\infty \leq \frac{1}{\mu_{\min}} \int_V |u_k-u|^2 d\mu,
\end{eqnarray*}
for any $2<r <\infty$, we have
\begin{eqnarray*}
\label{no014}
\int_V |u_k-u|^r d\mu \leq \|u_k-u\|^{r-2}_\infty \int_V |u_k-u|^2 d\mu \rightarrow 0 \;\;\text{as}\;\; k\rightarrow \infty.
\end{eqnarray*}
Therefore, up to a subsequence, $u_k \rightarrow u$ in $L^s(V)$ for all $2\leq r \leq \infty$. The proof is completed.
\qed
\vskip2mm
\noindent
\par
Let $E:=E_1 \times E_2$. Then $E$ is a Hilbert space with the inner product $\langle(u,v),(\phi,\psi)\rangle=(u,\phi)_{E_1}+(v,\psi)_{E_2}$. The corresponding norm is
$\|(u,v)\|_E=\left(\|u\|_{E_1}^2+\|v\|_{E_2}^2\right)^{\frac{1}{2}}$. The equivalent norm is $\|(u,v)\|=\|u\|_{E_1}+\|v\|_{E_2}$ which will be used in the sequel.  We define the functional $\Phi:E \rightarrow \R$ by
\begin{eqnarray}
\label{205}
\Phi(u,v)&=& \frac{1}{2}\int_V (|\Delta u|^2 + a_1|\nabla u|^2 +V_1(x)u^2)d\mu +\frac{1}{2}\int_V (|\Delta v|^2 + a_2|\nabla v|^2 +V_2(x)v^2)d\mu  \nonumber\\
&&+\frac{b_1}{4}\left(\int_V |\nabla u|^2d\mu\right)^2+\frac{b_2}{4}\left(\int_V |\nabla v|^2d\mu\right)^2- \int_V F(x,u,v)d\mu\nonumber\\
&=&\frac{1}{2}\|u\|_{E_1}^2+\frac{1}{2}\|v\|_{E_2}^2+\frac{b_1}{4}\left(\int_V |\nabla u|^2d\mu\right)^2+\frac{b_2}{4}\left(\int_V |\nabla v|^2d\mu\right)^2 - \int_V F(x,u,v)d\mu, \;\; (u,v) \in E.
\end{eqnarray}
\vskip2mm
\noindent
{\bf Lemma 2.2.} ({\cite{han2020, yang2024}}) {\it If $(V)$ and $(F_1)$ hold, then $\Phi \in C^1(E,\R)$ and
\begin{eqnarray}
\label{206}
\langle \Phi'(u,v),(\phi_1,\phi_2) \rangle
&=&(u,\phi_1)+ (v,\phi_2)
+ b_1 \int_V|\nabla u|^2d\mu \int_V \nabla u \nabla \phi_1 d\mu
+ b_2 \int_V|\nabla v|^2d\mu \int_V \nabla v \nabla \phi_2 d\mu \nonumber\\
&-& \int_V F_u(x,u,v)\phi_1 d\mu - \int_V F_v(x,u,v)\phi_2 d\mu,
\end{eqnarray}
for all $(u,v),\;(\phi_1,\phi_2) \in E$. Moreover, $\Psi': E\rightarrow E^*$ is compact, where $\Psi(u,v)= \int_V F(x,u,v)d\mu$.}

\vskip2mm
\par
By (\ref{206}), it is easy to obtain that
\begin{eqnarray}
\label{304}
\langle \Phi'(u,v),(u,v) \rangle
&=&\|u\|_{E_1}^2+ \|v\|_{E_2}^2
+ b_1 \left(\int_V |\nabla u|^2d\mu\right)^2 + b_2 \left(\int_V |\nabla v|^2d\mu\right)^2 \nonumber\\
&&- \int_V F_u(x,u,v)u d\mu - \int_V F_v(x,u,v)v d\mu, \;\; \forall (u,v) \in E.
\end{eqnarray}
\vskip2mm
\noindent
\par
Let $X$ be a Banach space. We say that $\Phi \in C^1(X,\R)$ satisfies $(C)$-condition if any sequence $\{u_n\}\subset X$ such that
\begin{eqnarray}
\label{208}
\Phi(u_n) \mbox{ is bounded}, \ \ \|\Phi'(u_n)\|_{X^*}(1+\|u_n\|)\rightarrow 0
\end{eqnarray}
has a convergent subsequence, where $\|\cdot\|_{X^*}$ is the norm in the dual space $X^*$.
\vskip2mm
\noindent
{\bf Lemma 2.3.} (\cite{Bartolo1983, Rabinowitz1986}) {\it Let $X$ be an infinite dimensional Banach space, $X=Y\bigoplus Z$, where $Y$ is finitely dimensional. If $\Phi \in C^1(X,\R)$ satisfies $(C)$-condition, and \\
$(I_1)$ $\Phi(0)=0$, $\Phi(-u)=\Phi(u)$ for all $u \in X$;\\
$(I_2)$ there exist constants $\rho, \alpha > 0$ such that $\Phi|_{\partial B_\rho \cap Z} \geq \alpha$;\\
$(I_3)$ for any finite dimensional subspace $\widetilde{X} \subset X$, there exists $R=R(\widetilde{X})>0$ such that $\Phi(u)\leq 0$ on $\widetilde{X}\setminus B_R$.\\
then $\Phi $ possesses an unbounded sequence of critical values.}

\section{Proofs for the  biharmonic-Kirchhoff case $(b_i>0, i=1,2)$}\label{section 3}
In this section, we aim to investigate system (\ref{101}) with $(b_i>0, i=1,2)$ and prove Theorem 1.1 and Theorem 1.2.  It is easy to see that $(u,v)\in E$ is a critical point of $\Phi$ if and only if
\begin{eqnarray*}
\label{no020}
(u,\phi_1)+ b_1 \int_V|\nabla u|^2 d\mu\int_V \nabla u \nabla \phi_1 d\mu- \int_V F_u(x,u,v)\phi_1 d\mu =0,
\end{eqnarray*}
and
\begin{eqnarray*}
\label{no021}
(v,\phi_2)+ b_2 \int_V|\nabla v|^2 d\mu\int_V \nabla v \nabla \phi_2 d\mu - \int_V F_v(x,u,v)\phi_2 d\mu=0,
\end{eqnarray*}
for all $(\phi_1,\phi_2)\in E$.
\par
In what follows, we will give some lemmas to prove our main theorems.
\vskip2mm
\noindent
{\bf Lemma 3.1.} {\it Under $(V)$, $(F_1)$, $(F_2)$, $(F_3)$ and  $(F_4)$,   any sequence $\{(u_n,v_n)\} \subset E$ satisfying
\begin{eqnarray}
\label{301}
\Phi(u_n,v_n)\mbox{ is bounded },\ \  \| \Phi'(u_n,v_n)\|(1+\|(u_n,v_n)\|) \rightarrow 0
\end{eqnarray}
is bounded in $E$.}
\vskip2mm
\noindent
{\bf Proof.} To prove the boundness of $\{(u_n,v_n)\}$, arguing by contradiction, assume that $\|(u_n,v_n)\|\rightarrow \infty$. Let $\phi_{1n}= \frac{u_n}{\|(u_n,v_n)\|}$ and $\phi_{2n}= \frac{v_n}{\|(u_n,v_n)\|}$. Then $\|\phi_{1n}\|_{E_1}\leq 1$, $\|\phi_{2n}\|_{E_2}\leq 1$, $\|\phi_{1n}\|_r\leq \gamma_{r,1}\|\phi_{1n}\|_{E_1}\le \gamma_{r,1}$ and $\|\phi_{2n}\|_r\leq \gamma_{r,2}\|\phi_{2n}\|_{E_2}\le \gamma_{r,2}$ for all $ r\geq2$.  Then, passing to a subsequence, we may assume that $\phi_{1n}\rightharpoonup \phi_1^*$ in $E_1$, $\phi_{2n}\rightharpoonup \phi_2^*$ in $E_2$. Furthermore, by Lemma 2.1, $\phi_{1n}\rightarrow \phi_1^* $, $\phi_{2n}\rightarrow \phi_2^* $ in $L^r(V)$ for all $r\geq2$, and $\phi_{1n}(x)\rightarrow \phi_1^*(x) $, $\phi_{2n}(x)\rightarrow \phi_2^*(x) $ for all $x \in V$.
\par
Firstly, we claim that $\phi^*_1 \not\equiv 0$ or $\phi^*_2\not\equiv 0$. Indeed, if $\phi_1^*\equiv \phi_2^* \equiv0$ for all $x\in V$, then $\phi_{1n},\;\phi_{2n}\rightarrow 0$ in $L^r(V)$ for all $r\geq2$, and $\phi_{1n}(x),\;\phi_{2n}(x)\rightarrow0$ for all $x\in V$.
Combining (\ref{205}) with (\ref{301}), we get
\begin{eqnarray}
\label{305}
2\limsup \limits_{n\rightarrow \infty} \int_V \frac{|F(x,u_n,v_n)|}{\|(u_n,v_n)\|^2} d\mu\ge \limsup \limits_{n\rightarrow \infty} \int_V \frac{|F(x,u_n,v_n)|}{\|u_n\|_{E_1}^2+\|v_n\|_{E_2}^2} d\mu\ge \frac{1}{2}.
\end{eqnarray}
For the constants $0\le d<e$, we denote
\begin{eqnarray}
\label{306}
\Omega_n[d,e)=\{x\in V: d\le |(u_n(x),v_n(x))|<e\}. \nonumber
\end{eqnarray}
Then, it follows from $(F_1)$ that
\begin{eqnarray}
\label{307}
\notag
&&\int_{\Omega_n[0,r_0)} \frac{|F(x,u_n,v_n)|}{\|(u_n,v_n)\|^2} d\mu \\
\notag
& \leq &\frac{1}{\|(u_n,v_n)\|^2} \int_{\Omega_n(0,r_0)} a(|(u_n,v_n)|)b(x) d\mu \\
&\leq &\frac{1}{\|(u_n,v_n)\|^2} \max\limits_{0\leq|(s,t)|\leq r_0}a(|(s,t)|) \int_V b(x)d\mu  \rightarrow 0, \mbox{ as } n\to \infty.
\end{eqnarray}
It follows from (\ref{205}), (\ref{206}) and (\ref{301}) that there exists a constant $c_1>0$ such that
\begin{eqnarray}
\label{302a}
 |\Phi(u_n,v_n)| \le c_1 ,\ \ \mbox{for all } n\in \mathbb N,
\end{eqnarray}
 and  for all large $n$,
\begin{eqnarray}
\label{302}
c_1 + 1 \geq \Phi(u_n,v_n) - \frac{1}{4}\langle\Phi'(u_n,v_n), (u_n,v_n)\rangle = \int_V \mathcal{F}(x,u_n,v_n)d\mu +\frac{1}{4}\|u_n\|^2+\frac{1}{4}\|v_n\|^2,
\end{eqnarray}
where $\mathcal{F}$ is defined in assumptions $(F_3)$. Set $k'=\frac{k}{k-1}$. Since $k> 1$, we get $2k'\in (2,\infty)$. Hence, from $(F_3)$ and (\ref{302}), we have
\begin{eqnarray}
\label{308}
\notag
&&\int_{\Omega_n[r_0,\infty)} \frac{|F(x,u_n,v_n)|}{\|(u_n,v_n)\|^2} d\mu \\
\notag
&=&\int_{\Omega_n[r_0,\infty)}\frac{F(x,u_n,v_n)}{|(u_n,v_n)|^2}|(\phi_{1n},\phi_{2n})|^2 d\mu \\
\notag
& \leq & \left(\int_{\Omega_n[r_0,\infty)}\left(\frac{F(x,u_n,v_n)}{|(u_n,v_n)|^2}\right)^k d\mu\right)^{\frac{1}{k}} \left(\int_{\Omega_n[r_0,\infty)}|(\phi_{1n},\phi_{2n})|^{2k'}d\mu\right)^{\frac{1}{k'}}\\
\notag
&\leq & c^{\frac{1}{k}}\left( \int_{\Omega_n[r_0,\infty)} \mathcal{F}(x,u_n,v_n)d\mu \right)^{\frac{1}{k}}\left(\int_{\Omega_n[r_0,\infty)}|(\phi_{1n},\phi_{2n})|^{2k'}d\mu\right)^{\frac{1}{k'}}\\
&\leq & [c(c_1+1)]^{\frac{1}{k}}\left(\int_{\Omega_n[r_0,\infty)}|(\phi_{1n},\phi_{2n})|^{2k'}d\mu\right)^{\frac{1}{k'}} \nonumber\\
&\leq & 2^{\frac{2k'-1}{k'}}[c(c_1+1)]^{\frac{1}{k}}\left(\int_{V}|\phi_{1n}|^{2k'}d\mu+\int_{V}|\phi_{2n}|^{2k'}d\mu\right)^{\frac{1}{k'}} \rightarrow 0, \mbox{ as } n\to \infty.
\end{eqnarray}
Combining $(F_0)$, (\ref{307}) with (\ref{308}), we get
\begin{eqnarray*}
\label{no023}
 \int_V \frac{|F(x,u_n,v_n)|}{\|(u_n,v_n)\|^2}d\mu
 = \int_{\Omega_n[0,r_0)} \frac{|F(x,u_n,v_n)|}{\|(u_n,v_n)\|^2}d\mu+ \int_{\Omega_n[r_0,\infty)} \frac{|F(x,u_n,v_n)|}{\|(u_n,v_n)\|^2}d\mu
 \rightarrow 0, \mbox{ as } n\to \infty,
\end{eqnarray*}
which contradicts with (\ref{305}). So $\phi^*_1 \not\equiv 0$ or $\phi^*_2\not\equiv 0$.
\par
Next, we prove $\{(u_n,v_n)\}$ is bounded in $E$. Without loss of generality, we assume that $\phi^*_1\not\equiv 0$. Set
$$A:= \{x \in V:\phi^*_1(x)\neq 0\}.$$
Then $A \not= \emptyset$. Thus, for all $x \in A$, we have
$$\lim_{n \rightarrow \infty } |u_n(x)|= +\infty.$$
Hence, $$A \subset \Omega_n[r_0,\infty) \;\;\text{for large}\;\;n \in \mathbb{N}.$$
It follows from (\ref{205}), (\ref{302a}), $(F_2)$ and Fatou's Lemma that
\begin{eqnarray}
\label{309}
\notag
0 &=& \lim_{n\rightarrow\infty}\frac{-c}{\|(u_n,v_n)\|^4}\\
\notag
& \le & \lim_{n\rightarrow\infty}\frac{\Phi(u_n,v_n)}{\|(u_n,v_n)\|^4}\\
\notag
&\leq & \lim_{n\rightarrow \infty} \left[\max\left\{\frac{b_1}{4a_1^2},\frac{b_2}{4a_2^2}\right\}-\int_V \frac{F(x,u_n,v_n)}{\|(u_n,v_n)\|^4} d\mu \right]\\
\notag
&=& \lim_{n\rightarrow \infty} \left[\max\left\{\frac{b_1}{4a_1^2},\frac{b_2}{4a_2^2}\right\}- \int_{\Omega_n[0,r_0)} \frac{F(x,u_n,v_n)}{\|(u_n,v_n)\|^4} d\mu-\int_{\Omega_n[r_0,\infty)} \frac{F(x,u_n,v_n)}{|(u_n,v_n)|^4} |(\phi_{1n},\phi_{2n})|^4 d\mu \right]\\
\notag
&\leq & \lim\sup_{n\rightarrow\infty}\left[\max\left\{\frac{b_1}{4a_1^2},\frac{b_2}{4a_2^2}\right\}+\frac{1}{\|(u_n,v_n)\|^4}\max\limits_{0\leq|(s,t)|\leq r_0}a(|(s,t)|) \int_V b(x)d\mu\right.\\
\notag
&& \left.-\int_{\Omega_n[r_0,\infty)} \frac{F(x,u_n,v_n)}{|(u_n,v_n)|^4} |(\phi_{1n},\phi_{2n})|^4d\mu\right]\\
\notag
&=& \max\left\{\frac{b_1}{4a_1^2},\frac{b_2}{4a_2^2}\right\}-\liminf_{n\rightarrow \infty}\int_{\Omega_n[r_0,\infty)} \frac{F(x,u_n,v_n)}{|(u_n,v_n)|^4} |(\phi_{1n},\phi_{2n})|^4d\mu \\
\notag
&\le&  \max\left\{\frac{b_1}{4a_1^2},\frac{b_2}{4a_2^2}\right\} -\int_{V} \liminf_{n\rightarrow \infty}\frac{F(x,u_n,v_n)}{|(u_n,v_n)|^4} \chi_{\Omega_n[r_0,\infty)}|(\phi_{1n},\phi_{2n})|^4d\mu \\
\notag
&\leq & \max\left\{\frac{b_1}{4a_1^2},\frac{b_2}{4a_2^2}\right\}-\int_{A} \liminf_{n\rightarrow \infty} \frac{F(x,u_n,v_n)}{|(u_n,v_n)|^4} \chi_{\Omega_n[r_0,\infty)}|(\phi_{1n},\phi_{2n})|^4d\mu \\
&=& -\infty,
\end{eqnarray}
where $\chi_{\Omega_n[r_0,\infty)}=\begin{cases}1,\ \ x\in \Omega_n[r_0,\infty)\\
0,\ \ x\in V \backslash \Omega_n[r_0,\infty) \end{cases}$, which is a contradiction. Thus $\{(u_n,v_n)\}$ is bounded in $E$.
\vskip2mm
\noindent
{\bf Lemma 3.2.} {\it Under $(V)$, $(F_1)$, $(F_2)$, $(F_3)$ and  $(F_4)$, any sequence $\{(u_n,v_n)\} \subset E$ satisfying (\ref{301}) has a convergent subsequence in $E$.}
\vskip2mm
\noindent
{\bf Proof.} From Lemma 3.1, we know that $\{(u_n,v_n)\}$ is bounded in $E$, that is, there exists a constant $M_1$ such that $\|(u_n,v_n)\|_{E_1} \leq M_1$. Then, by Lemma 2.1, $\|u_n\|_\infty \leq M_1 \gamma_{\infty,1}$ and $\|v_n\|_\infty \leq M_1 \gamma_{\infty,2}$. Going if necessary to a subsequence, we can assume that $(u_n,v_n) \rightharpoonup (u,v)$ in $E$. By Lemma 2.1, $u_n\rightarrow u$ and $v_n\rightarrow v$ in $L^s(V)$ for all $2 \leq s \leq +\infty$. Thus
\begin{eqnarray}
\label{310}
\notag
&&\int_V |F_u(x,u_n,v_n)-F_u(x,u,v)||u_n-u|d\mu\\
\notag
& \leq & \int_V (|F_u(x,u_n,v_n)|+|F_u(x,u,v)|)|u_n-u|d\mu\\
\notag
& \leq & \int_V 2 a(|(u_n,v_n)|)b(x)|u_n-u|d\mu\\
\notag
& \leq & 2 \max\limits_{|(s,t)|\leq M_1\gamma_{\infty,1}+ M_1\gamma_{\infty,2}} a(|(s,t)|) \int_V b(x)d\mu \|u_n-u\|_\infty\\
&\rightarrow & 0,\;\; \text{as}\;\; n \rightarrow\infty.
\end{eqnarray}
Similarly, we also obtain
\begin{eqnarray}
\label{320}
\int_V |F_v(x,u_n,v_n)-F_v(x,u,v)||v_n-v|d\mu \rightarrow  0,\;\; \text{as}\;\; n \rightarrow\infty.
\end{eqnarray}
Moreover, by Lemma 2.1, we have $u_n(x)\to u(x)$ and $v_n(x)\to v(x)$ for all $x\in V$. Then by  the definition of $\nabla u$, it is easy to see that
$\nabla u_n(x)\to u_n(x)$ for all $x\in V$, as $n\to \infty$. Then by the boundness of $\{u_n\}$ in $E$, we have
\begin{eqnarray}
\label{311a}
&      &   b_1\left( \int_V|\nabla u|^2 d\mu \int_V \nabla u \nabla (u_n-u) d\mu- \int_V|\nabla u_n|^2 d\mu \int_V \nabla u_n \nabla (u_n-u) d\mu \right)\nonumber\\
& \le  &   b_1\left( \int_V|\nabla u|^2 d\mu \int_V |\nabla u| |\nabla (u_n-u)| d\mu+\int_V|\nabla u_n|^2 d\mu \int_V |\nabla u_n| |\nabla (u_n-u)| d\mu \right)\nonumber\\
& \le  &   b_1\left[ \left(\int_V|\nabla u|^2 d\mu\right)^{\frac{3}{2}} +\left(\int_V|\nabla u_n|^2 d\mu\right)^{\frac{3}{2}} \right]\left(\int_V |\nabla u_n-\nabla u|^2 d\mu\right)^{\frac{1}{2}}\nonumber\\
&\to & 0,\ \ \mbox{as } n\to \infty .
\end{eqnarray}
Observe that
\begin{eqnarray}
\label{311}
\notag
&&\|u_n-u\|^2_{E_1}+\|v_n-v\|^2_{E_2} \\
\notag
&=&  \langle \Phi'(u_n,v_n)-\Phi'(u,v),(u_n-u,v_n-v)\rangle \\
\notag
&&+ \int_V [F_u(x,u_n,v_n)-F_u(x,u,v)](u_n-u) d\mu + \int_V [F_v(x,u_n,v_n)-F_v(x,u,v)](v_n-v) d\mu \\
\notag
&&+ b_1\left( \int_V|\nabla u|^2 d\mu \int_V \nabla u \nabla (u_n-u) d\mu- \int_V|\nabla u_n|^2 d\mu \int_V \nabla u_n \nabla (u_n-u) d\mu \right)\\
&&+ b_2 \left( \int_V|\nabla v|^2 d\mu \int_V \nabla v \nabla (v_n-v) d\mu-  \int_V|\nabla v_n|^2 d\mu \int_V \nabla v_n \nabla (v_n-v) d\mu \right).
\end{eqnarray}
By (\ref{301}), it is clear that $\langle \Phi'(u_n,v_n)-\Phi'(u,v),(u_n-u,v_n-v)\rangle \to 0$ as $n\rightarrow \infty $. From (\ref{310}), (\ref{320}), (\ref{311a}) and (\ref{311}), we get
\begin{eqnarray*}
\label{no027}
\|u_n-u\|_{E_1} \rightarrow 0 \;\;\text{as}\;\; n\rightarrow\infty.
\end{eqnarray*}
and
\begin{eqnarray*}
\label{no028}
\|v_n-v\|_{E_2} \rightarrow 0 \;\;\text{as}\;\; n\rightarrow\infty.
\end{eqnarray*}
The proof is completed.\qed

\vskip2mm
\noindent
{\bf Lemma 3.3.} {\it Under $(V)$, $(F_1)$, $(F_2)$, $(F_3)$ and  $(F_5)$, any sequence $\{(u_n,v_n)\} \subset E$ satisfying (\ref{301}) has a convergent subsequence in $E$.}
\vskip2mm
\noindent
{\bf Proof.} Firstly, we prove that $\{(u_n,v_n)\}$ is bounded in $E$. To prove the boundedness of $\{(u_n,v_n)\}$, arguing by contradiction, suppose that $\|(u_n,v_n)\|\rightarrow\infty$. By (\ref{301}), (\ref{304}) and $(F_4)$, for large $n$, there exists a constant $c>0$ such that
\begin{eqnarray}
\label{312}
\notag
c+1 &\geq & \Phi(u_n,v_n)-\frac{1}{\mu}\langle\Phi'(u_n,v_n),(u_n,v_n)\rangle\\
\notag
&=& \frac{\mu-2}{2\mu}(\|u_n\|^2_{E_1}+\|v_n\|^2_{E_2})+\frac{(\mu-4)b_1}{4\mu}\left(\int_V |\nabla u_n|^2d\mu\right)^2 + \frac{(\mu-4)b_2}{4\mu}\left(\int_V |\nabla v_n|^2d\mu\right)^2\\
\notag
&&+ \int_V\left[\frac{1}{\mu}F_u(x,u_n,v_n)u_n+\frac{1}{\mu}F_v(x,u_n,v_n)v_n-F(x,u_n,v_n)\right]d\mu\\
\notag
&\geq & \frac{\mu-2}{2\mu}(\|u_n\|^2_{E_1}+\|v_n\|^2_{E_2})+\frac{(\mu-4)b_1}{4\mu}\left(\int_V |\nabla u_n|^2d\mu\right)^2 + \frac{(\mu-4)b_2}{4\mu}\left(\int_V |\nabla v_n|^2d\mu\right)^2 \\
\notag
&&-\frac{\sigma}{\mu}(\|u_n\|^2_2+\|v_n\|^2_2)\\
\notag
&\geq & \frac{\mu-2}{2\mu}(\|u_n\|^2_{E_1}+\|v_n\|^2_{E_2})- \frac{\sigma}{\mu}(\gamma^2_{2,1}\|u_n\|^2_{E_1}+\gamma^2_{2,2}\|v_n\|^2_{E_2})\\
& \geq & \left(\frac{\mu-2}{2\mu}-  \frac{\sigma}{\mu} \max\{\gamma^2_{2,1},\gamma^2_{2,2}\}\right)(\|u_n\|^2_{E_1}+\|v_n\|^2_{E_2})
\end{eqnarray}
Note that $\sigma\in \left(0, \frac{\mu-2}{2\max\{\gamma_{2,1}^2,\gamma_{2,2}^2\}}\right)$. Thus, (\ref{312}) implies that $\{(u_n,v_n)\}$ is bounded in $E$. The rest proof is the same as that in Lemma 3.2.
\qed
\vskip2mm
\noindent
{\bf Lemma 3.4.} {\it Under $(V)$, $(F_1)$ and $(F_2)$, for any finite dimensional subspace $\widetilde{E}= \widetilde{E_1}\times \widetilde{E_2}\subset E$, where $\widetilde{E_1}\subset E_1$ and $\widetilde{E_2}\subset E_2$, there holds
\begin{eqnarray}
\label{314}
\Phi (u,v) \rightarrow -\infty,\; \text{as}\; \|(u,v)\| \rightarrow \infty,\; \forall (u,v) \in \widetilde{E}.
\end{eqnarray}}
\noindent
{\bf Proof.} Arguing indirectly, assume that for some sequence $\{(u_n,v_n)\} \subset \widetilde{E}$ with $\|(u_n,v_n)\| \rightarrow \infty$, there exists $M_1 > 0$ such that $\Phi (u_n,v_n) \geq -M_1$ for all $n \in \mathbb{N}$. Next, we argue from  three cases.
 \par
 {\bf Case (I): assume that $\|u_n\|_{E_1}\to \infty$ and  $\|v_n\|_{E_2}\to \infty$.} Let $\varphi_{1n} = \frac{u_n}{\|u_n\|_{E_1}}$ and $\varphi_{2n} = \frac{v_n}{\|v_n\|_{E_2}}$. Then $\|\varphi_{1n}\|_{E_1} = 1$ and $\|\varphi_{2n}\|_{E_2} =1$. Then we may assume that there exists a subsequence, still denoted by $\{\varphi_{1n}\}$ and $\{\varphi_{2n}\}$, such that $\varphi_{1n} \rightharpoonup \varphi_1^*$ and $\varphi_{2n} \rightharpoonup \varphi_2^*$ in $E$. Since $\widetilde{E}$ is finite dimensional, then $\varphi_{1n} \rightarrow \varphi_1^*$, $\varphi_{2n} \rightarrow \varphi_2^* \in \widetilde{E}$ in $E$, and $\varphi_{1n}(x) \rightarrow \varphi_1^*(x)$, $\varphi_{2n}(x) \rightarrow \varphi_2^*(x)$ for all $x \in V$, and so $\|\varphi_1^*\|_{E_1} = 1$ and $\|\varphi_2^*\|_{E_2} = 1$.
 Set
$$
B:= \{x \in V:\varphi^*_1(x)\neq 0\}.
$$
 We get $B \not\equiv \emptyset$. For all $x \in B$, we have
$$\lim_{n \rightarrow \infty } |u_n(x)|= +\infty.$$
Hence, $$B \subset \Omega_n[r_0,\infty) \;\;\text{for large}\;\;n \in \mathbb{N}.$$
Similar to the argument of (\ref{309}), we have the following contradiction:
\begin{eqnarray}
\label{a1}
\notag
0&=&\lim_{n\rightarrow\infty}\frac{-M_1}{\|(u_n,v_n)\|^4} \\
\notag
& \le & \lim_{n\rightarrow\infty}\frac{\Phi(u_n,v_n)}{\|(u_n,v_n)\|^4}\\
\notag
&\leq & \lim_{n\rightarrow \infty} \left[\max\left\{\frac{b_1}{4a_1^2},\frac{b_2}{4a_2^2}\right\}-\int_V \frac{F(x,u_n,v_n)}{\|(u_n,v_n)\|^4} d\mu \right]\\
\notag
&\leq &  \max\left\{\frac{b_1}{4a_1^2},\frac{b_2}{4a_2^2}\right\}-\liminf_{n\rightarrow \infty}\int_{\Omega_n[r_0,\infty)} \frac{F(x,u_n,v_n)}{|(u_n,v_n)|^4} \frac{|(\varphi_{1n}\|u_n\|_{E_1},\varphi_{2n}\|v_n\|_{E_2})|^4}{\|(u_n,v_n)\|^4}d\mu \\
\notag
&\leq &  \max\left\{\frac{b_1}{4a_1^2},\frac{b_2}{4a_2^2}\right\}-\frac{1}{8}\liminf_{n\rightarrow \infty}\int_{\Omega_n[r_0,\infty)} \frac{F(x,u_n,v_n)}{|(u_n,v_n)|^4} \min\{|\varphi_{1n}|^4,|\varphi_{2n}|^4\}d\mu \\
\notag
&\le&  \max\left\{\frac{b_1}{4a_1^2},\frac{b_2}{4a_2^2}\right\} -\frac{1}{8}\int_{V} \liminf_{n\rightarrow \infty}\frac{F(x,u_n,v_n)}{|(u_n,v_n)|^4} \chi_{\Omega_n[r_0,\infty)}\min\{|\varphi_{1n}|^4,|\varphi_{2n}|^4\}d\mu \\
\notag
&\leq & \max\left\{\frac{b_1}{4a_1^2},\frac{b_2}{4a_2^2}\right\}-\frac{1}{8}\int_{B} \liminf_{n\rightarrow \infty} \frac{F(x,u_n,v_n)}{|(u_n,v_n)|^4} \chi_{\Omega_n[r_0,\infty)}\min\{|\varphi_{1n}|^4,|\varphi_{2n}|^4\}d\mu \\
&=& -\infty,
\end{eqnarray}
where we have used the inequality that
$$
|(\varphi_{1n}\|u_n\|_{E_1},\varphi_{2n}\|v_n\|_{E_2})|^4 \ge |\varphi_{1n}|^4 \|u_n\|_{E_1}^4+  |\varphi_{2n}|^4 \|v_n\|_{E_2}^4 \ge \frac{1}{8}\min\{|\varphi_{1n}|^4,|\varphi_{2n}|^4\}\|(u_n,v_n)\|^4.
$$
 \par
  {\bf Case (II): assume that $\|u_n\|_{E_1}\to \infty$ and  $\|v_n\|_{E_2}$ is bounded. } Let $\hat{\varphi}_{1n} = \frac{u_n}{\|u_n\|_{E_1}}$. Then $\|\hat{\varphi}_{1n}\|_{E_1} = 1$ and there exists a constant $M_2>0$ such that $\|v_n\|_{E_2} \le M_2$. Then we may assume that there exists a subsequence, still denoted by $\{\hat{\varphi}_{1n}\}$, such that $\hat{\varphi}_{1n} \rightharpoonup \hat{\varphi}_1^*$ in $E_1$. Since $\widetilde{E_1}$ is finite dimensional, then $\hat{\varphi}_{1n} \rightarrow \hat{\varphi}_1^* \in \widetilde{E_1}$ in $E$, and $\hat{\varphi}_{1n}(x) \rightarrow \hat{\varphi}_1^*(x)$ for all $x \in V$, and so $\|\hat{\varphi}_1^*\|_{E_1} = 1$.
 Set
$$
D:= \{x \in V:\hat{\varphi}^*_1(x)\neq 0\}.
$$
 We get $D \not\equiv \emptyset$. For all $x \in D$, we have
$$\lim_{n \rightarrow \infty } |u_n(x)|= +\infty.$$
Hence, $$D \subset \Omega_n[r_0,\infty) \;\;\text{for large}\;\;n \in \mathbb{N}.$$
Similar to the argument of (\ref{309}),  we have the following contradiction:
\begin{eqnarray}
\label{b4}
\notag
0&=&\lim_{n\rightarrow\infty}\frac{-M_1}{\|(u_n,v_n)\|^4} \\
\notag
& \le & \lim_{n\rightarrow\infty}\frac{\Phi(u_n,v_n)}{\|(u_n,v_n)\|^4}\\
\notag
&\leq & \lim_{n\rightarrow \infty} \left[\max\left\{\frac{b_1}{4a_1^2},\frac{b_2}{4a_2^2}\right\}-\int_V \frac{F(x,u_n,v_n)}{\|(u_n,v_n)\|^4} d\mu \right]\\
\notag
&\leq &  \max\left\{\frac{b_1}{4a_1^2},\frac{b_2}{4a_2^2}\right\}-\liminf_{n\rightarrow \infty}\int_{\Omega_n[r_0,\infty)} \frac{F(x,u_n,v_n)}{|(u_n,v_n)|^4} \frac{|(\hat{\varphi}_{1n}\|u_n\|_{E_1},|v_n|)|^4}{\|(u_n,v_n)\|^4}d\mu \\
\notag
&\leq &  \max\left\{\frac{b_1}{4a_1^2},\frac{b_2}{4a_2^2}\right\}-\frac{1}{16}\liminf_{n\rightarrow \infty}\int_{\Omega_n[r_0,\infty)} \frac{F(x,u_n,v_n)}{|(u_n,v_n)|^4} |\hat{\varphi}_{1n}|^4d\mu \\
\notag
&\le&  \max\left\{\frac{b_1}{4a_1^2},\frac{b_2}{4a_2^2}\right\} -\frac{1}{16}\int_{V} \liminf_{n\rightarrow \infty}\frac{F(x,u_n,v_n)}{|(u_n,v_n)|^4} \chi_{\Omega_n[r_0,\infty)}|\hat{\varphi}_{1n}|^4d\mu \\
\notag
&\leq & \max\left\{\frac{b_1}{4a_1^2},\frac{b_2}{4a_2^2}\right\}-\frac{1}{16}\int_{D} \liminf_{n\rightarrow \infty} \frac{F(x,u_n,v_n)}{|(u_n,v_n)|^4} \chi_{\Omega_n[r_0,\infty)}|\hat{\varphi}_{1n}|^4d\mu \\
&=& -\infty,
\end{eqnarray}
where we have used the inequality that
$$
\frac{|(\hat{\varphi}_{1n}\|u_n\|_{E_1},|v_n|)|^4}{(\|u_n\|_{E_1}+\|v_n\|_{E_1})^4} \ge\frac{\hat{\varphi}_{1n}^4\|u_n\|_{E_1}^4}{(\|u_n\|_{E_1}+M_2)^4} \ge\frac{\hat{\varphi}_{1n}^4\|u_n\|_{E_1}^4}{(\|u_n\|_{E_1}+\|u_n\|_{E_1})^4} =\frac{\hat{\varphi}_{1n}^4}{2^4} \mbox{ for all sufficiently large }n.
$$

\vskip2mm
\par
{\bf Case (III): assume that $\|u_n\|_{E_1}$ is bounded and  $\|v_n\|_{E_2}\to \infty$.} This case is similar to Case (II) by exchanging $u_n$ with $v_n$. We omit the details. The proof is completed.
\qed

\vskip2mm
\noindent
{\bf Corollary 3.5.} {\it Under $(V)$, $(F_1)$ and $(F_2)$, for any finite dimensional subspace $\widetilde{E} \subset E$, there exists $R=R(\widetilde{E})>0$, such that}
$
\Phi (u,v) \leq 0,\; \forall (u,v) \in \widetilde{E},\; \|(u,v)\| \geq R.
$

\vskip2mm
\par
Let $\{e_{ij}\}_{j=1}^\infty, i=1,2$ are two total orthonormal  bases of $E_i, i=1,2$, respectively. For any given $k\in \mathbb N$, let
\begin{eqnarray*}
\label{316}
X_j = (\R e_{1j}, \R e_{2j}),\; Y_k = \bigoplus\limits_{j=1}^k X_j = \left(\bigoplus\limits_{j=1}^k \R e_{1j}, \bigoplus\limits_{j=1}^k \R e_{2j}\right), \;Z_k = \bigoplus\limits_{j=k+1}^\infty X_j =\left(\bigoplus\limits_{j=k+1}^\infty \R e_{1j}, \bigoplus\limits_{j=k+1}^\infty \R e_{2j}\right) ,\; k \in \mathbb{Z}.
\end{eqnarray*}
Then $E=Y_k\bigoplus Z_k$ and $Y_k$ is of finite dimension.
\vskip2mm
\noindent
{\bf Lemma 3.6.} {\it Under $(V)$, for $ r\geq2 $, it holds that
\begin{eqnarray}
\label{317}
\beta_k := \sup\limits_{(u,v) \in Z_k,\|(u,v)\|=1} (\|u\|_r+\|v\|_r) \rightarrow 0,\; \mbox{as } k \rightarrow \infty.
\end{eqnarray}}
\vskip2mm
\noindent
{\bf Proof.}  It is easy to see  that $0 < \beta_{k+1} \leq \beta_k$. Hence, $\beta_k \rightarrow \beta \geq 0 \;(k \rightarrow \infty)$. By the definition of $\beta_k$, we obtain that  for every $k \in \mathbb{N}$, there exists $(u_k,v_k) \in Z_k$ with $\|(u_k,v_k)\| = 1$ such that
\begin{eqnarray}\label{a3}
\|u\|_r+\|v\|_r> \frac{\beta_k}{2}
\end{eqnarray}
 For any given $(\phi_1,\phi_2) \in E$, writing $\phi_1 = \Sigma_{j=1} ^\infty c_{1j} e_{1,j}$ and $\phi_2 = \Sigma_{j=1} ^\infty c_{2j} e_{2,j}$, where $c_{1j}=(\phi_1, e_{1,j})_{E_1}$ and
$c_{2j}=(\phi_2, e_{2,j})_{E_2}$. By Parseval equality,  we have $\|\phi_1\|_{E_1}^2=\sum_{i=1}^\infty c_{1j}^2$ and $\|\phi_2\|_{E_2}^2=\sum_{i=1}^\infty c_{2j}^2$. Furthermore, by $\phi_1 \in E_1$ and $\phi_2 \in E_2$, we get
\begin{eqnarray}\label{a2}
 \left(\sum_{i=k+1}^\infty c_{1i}^2\right)^{\frac{1}{2}} + \left(\sum_{j=k+1}^\infty c_{2j}^2\right)^{\frac{1}{2}}\rightarrow 0, \mbox{ as  } k\rightarrow \infty.
\end{eqnarray}
Then, it follows from  the Cauchy-Schwartz inequality and (\ref{a2}) that
\begin{eqnarray*}
\label{no028}
|\langle (u_k,v_k),(\phi_1,\phi_2)\rangle|&=&\left|(u_k,\sum_{j=1}^\infty c_{1j} e_{1,j})\right|+ \left|(v_k,\sum_{j=1}^\infty c_{2j} e_{2,j})\right|\\
&=&\left|(u_k,\sum_{i=k+1}^\infty c_{1j} e_{1,j})\right|+ \left|(v_k,\sum_{j=k+1}^\infty c_{2j} e_{2,j})\right|\\
\notag
&\leq& \|u_k\|_{E_1} \cdot \left\|\sum_{i=k+1}^\infty c_{1j} e_{1,j}\right\|_{E_1}+ \|v_k\|_{E_2} \cdot \left\|\sum_{j=k+1}^\infty c_{2j} e_{2,j}\right\|_{E_2}\\
 &=& \left(\sum_{i=k+1}^\infty c_{1j}^2\right)^{\frac{1}{2}} + \left(\sum_{j=k+1}^\infty c_{2j}^2\right)^{\frac{1}{2}}\to 0\ \ \mbox{as } k\to \infty.
\end{eqnarray*}
 which, together with the Riesz lemma, implies that $(u_k,v_k)\rightharpoonup (0,0)$ in $E$ and then $u_k\rightharpoonup 0$ in $E_1$ and $v_k\rightharpoonup 0$ in $E_2$. By Lemma 2.1, the compact embedding of $E_i\hookrightarrow L^r(V)$ $(r\geq2, i=1,2)$ implies that $u_k\rightarrow 0$ and $v_k\rightarrow 0$ in $L^r(V)$. Hence, letting $k\rightarrow\infty$, it follows from (\ref{a3}) that $\beta=0$. Thus the proof is completed.
\qed

\vskip2mm
\noindent
{\bf Lemma 3.7.} {\it  Under $(V)$, $(F_1)$ and $(C_1)$, for any given large $k$, there exist constants $\rho,\alpha > 0$ such that $\Phi |_{\partial B_\rho \cap Z_k \geq \alpha}$.}
\vskip2mm
\noindent
{\bf Proof.}
It follows from $(C_1)$ that for any given $\varepsilon>0$,  there exists a $0<\delta<\frac{\max\{\gamma_{\infty, 1}, \gamma_{\infty,2}\}}{2^{\frac{p+1}{p-2}}(h+\varepsilon)\|b\|_\infty^{\frac{1}{p-2}}}$ such that
 \begin{eqnarray}
\label{a4} a(|(s,t)|)\le (h+\varepsilon)|(t,s)|^2\le 2(h+\varepsilon)(|t|^2+|s|^2), \mbox{ for }\forall |(t,s)|\le \delta.
\end{eqnarray}
For all $(u,v)\in Z_k$ with $\|(u,v)\|=\rho=\frac{\delta}{\max\{\gamma_{\infty, 1}, \gamma_{\infty,2}\}}$, by Lemma 2.1, we have
 \begin{eqnarray} \label{a5}
  |(u(x),v(x))|\le \|u\|_\infty+\|v\|_\infty\le \gamma_{\infty, 1} \|u\|_{E_1}+\gamma_{\infty,2}\|v\|_{E_2}\le \max\{\gamma_{\infty, 1}, \gamma_{\infty,2}\} \|(u,v)\|\le \delta.
\end{eqnarray}
Moreover,  by Lemma 3.6,  for any given large $k$, it holds that
\begin{eqnarray}
\label{318}
\left(\|u\|_2+\|v\|_2\right)^2 \leq \frac{1}{16(h+\varepsilon)\|b\|_\infty} \|(u,v)\|^2,\; \forall (u,v) \in Z_k.
\end{eqnarray}
Note that the fact that $b\in L^1(V)$ implies that $b\in L^\infty(V)$. Then for all  $(u,v)\in Z_k$ with $\|(u,v)\|=\rho$, by $(F_1)$, (\ref{a4}), (\ref{a5}) and (\ref{318}), we obtain that
\begin{eqnarray*}
\label{no030}
\Phi(u,v) &=& \frac{1}{2}\|u\|^2_{E_1}+\frac{1}{2}\|v\|^2_{E_2}+\frac{b_1}{4}\left(\int_V |\nabla u|^2d\mu\right)^2+\frac{b_2}{4}\left(\int_V |\nabla v|^2d\mu\right)^2 - \int_V F(x,u,v)d\mu\\
\notag
&\geq & \frac{1}{2}\|u\|^2_{E_1}+\frac{1}{2}\|v\|^2_{E_2} - \int_V a(|(u,v)|)b(x)d\mu \\
\notag
&\geq & \frac{1}{2}\|u\|^2_{E_1}+\frac{1}{2}\|v\|^2_{E_2} - 2(h+\varepsilon)\|b\|_\infty (\|u\|^2_2+\|v\|_2^2)\\
\notag
&\geq & \frac{1}{2}\|u\|^2_{E_1}+\frac{1}{2}\|v\|^2_{E_2} - 2(h+\varepsilon)\|b\|_\infty (\|u\|_2+\|v\|_2)^2\\
\notag
&\geq & \frac{1}{2}\|u\|^2_{E_1}+\frac{1}{2}\|v\|^2_{E_2}-2(h+\varepsilon)\|b\|_\infty\cdot  \frac{1}{16(h+\varepsilon)\|b\|_\infty}  \|(u,v)\|^2\\
\notag
& \geq & \frac{1}{2}\cdot \frac{1}{2} \|(u,v)\|^2 -\frac{1}{8}\|(u,v)\|^2\\
\notag
&=& \frac{1}{8}\rho^2:=\alpha > 0.
\end{eqnarray*}
Thus, the proof is complete.
\qed
\vskip2mm
\noindent
{\bf Lemma 3.8.} {\it Under $(V)$, $(F_1)$ and $(C_2)$,  $\Phi |_{\partial B_{\rho_*} \cap Z_k} \geq \alpha_*$, where $\rho_*$ and $\alpha_*$ are given in $(C_2)$.}
\vskip2mm
\noindent
{\bf Proof.} By  $(F_1)$ and $(C_2)$, for $(u,v) \in Z_k$ with $\|(u,v)\|=\rho$, we have
\begin{eqnarray*}
\label{no030}
\Phi(u,v) &=& \frac{1}{2}\|u\|^2_{E_1}+\frac{1}{2}\|v\|^2_{E_2}+\frac{b_1}{4}\left(\int_V |\nabla u|^2d\mu\right)^2+\frac{b_2}{4}\left(\int_V |\nabla v|^2d\mu\right)^2 - \int_V F(x,u,v)d\mu\\
\notag
&\geq & \frac{1}{2}\|u\|^2_{E_1}+\frac{1}{2}\|v\|^2_{E_2}- \int_V a(|(u,v)|)b(x)d\mu \\
\notag
&\geq & \frac{1}{2}\|u\|^2_{E_1}+\frac{1}{2}\|v\|^2_{E_2} -  \max\limits_{|s|+|t|\leq(\gamma_{\infty,1}+\gamma_{\infty,2})\rho} a(|(s,t)|) \int_V b(x)d\mu \\
\notag
& \geq & \frac{1}{2}\cdot \frac{1}{2} (\|u\|_{E_1}+\|v\|_{E_2})^2 -\max\limits_{|s|+|t|\leq(\gamma_{\infty,1}+\gamma_{\infty,2})\rho} a(|(s,t)|) \int_V b(x)d\mu \\
\notag
&=& \alpha_* > 0.
\end{eqnarray*}
Thus, the proof is completed.
\qed
\vskip2mm
\noindent
{\bf Proof of Theorem 1.1.} Let $X=E,\;Y=Y_k,\;Z=Z_k$.  By Lemmas 3.1, 3.2, 3.7 and Corollary 3.5, all conditions of Lemma 2.4 are satisfied. Thus, problem (\ref{101}) possesses infinitely many nontrivial solutions $\{(u_k,v_k)\}$ such that (\ref{102}) holds.\qed
\vskip2mm
\noindent
{\bf Proof of Theorem 1.2.} Let $X=E,\;Y=Y_k,\;Z=Z_k$. By Lemmas 3.1, 3.2, 3.8 and Corollary 3.5, all conditions of Lemma 2.4 are satisfied. Thus, problem (\ref{101}) possesses infinitely many nontrivial solutions $\{(u_k,v_k)\}$ such that (\ref{102}) holds. \qed

\section{The results corresponding to the  biharmonic case  $(b_i=0, i=1,2)$}\label{section 4}
In this section, corresponding to Theorem 1.1 and Theorem 1.2, we present some results for system (\ref{101}) with $b_i=0, i=1,2$, that is,
 \begin{eqnarray}
\label{101a}
   \begin{cases}
 \Delta^2 u- a_1\Delta u + V_1(x)u=F_u(x,u,v),\;\;\;\;\hfill x\in V,\\
 \Delta^2 v- a_2\Delta v + V_2(x)v=F_v(x,u,v),\;\;\;\;\hfill x\in V,\\
   \end{cases}
\end{eqnarray}
 and their proofs are almost same as Theorem 1.1 and Theorem 1.2 by simply changing $4$ into $2$ in (\ref{302}), (\ref{309}), (\ref{312}), (\ref{a1}), (\ref{b4}) and so on. Hence, we omit the detail proofs.
Next, we state the results.
\vskip2mm
\noindent
{\bf Theorem 4.1.} {\it Suppose that $(V)$, $(F_0), (F_1), (C_1) (\mbox{or } (C_2)),  (F_3)$ and $(F_6)$ and the following conditions hold:\\
$(F_2)'$ $\lim_{|(s,t)|\rightarrow\infty} \frac{F(x,s,t)}{|(s,t)|^2}=+\infty$ uniformly for all $x \in V$;\\
$(F_4)'$ $\mathcal{F}(x,s,t)= \frac{1}{2}\left(F_s(x,s,t)s+F_t(x,s,t)t\right)-F(x,s,t)\geq 0$, and there exist $c_1>0$ and $k > 1$ such that
$$|F(x,s,t)|^k \leq c_1 |(s,t)|^{2k}\mathcal{F}(x,s,t)\;\text{for all}\;x \in V,\;(s,t) \in \R^2,\;|(s,t)| \geq r_0;$$
Then problem (\ref{101a}) has infinitely many solutions $\{(u_k,v_k)\}_{k=1}^\infty$ such that
\begin{eqnarray}
\label{102a}
&&\frac{1}{2} \int_V (|\Delta u_k|^2+ a_1|\nabla u_k|^2+ V_1(x)u_k^2)d\mu +\frac{1}{2} \int_V (|\Delta v_k|^2+ a_2|\nabla v_k|^2+ V_2(x)v_k^2)d\mu\nonumber\\
&&- \int_V F(x,u_k,v_k)d\mu \rightarrow + \infty, \;\;\text{as} \;k \rightarrow + \infty .
 \end{eqnarray}}
\noindent
{\bf Theorem 4.2.} {\it Suppose that $(V)$, $(F_0), (F_1), (C_1) (\mbox{or } (C_2)), (F_2)', (F_3)$ and $(F_6)$ and the following condition hold:\\
  $(F_5)'$ there exist $\mu > 2$ and $\sigma\in \left(0, \frac{\mu-2}{2\max\{\gamma_{2,1}^2,\gamma_{2,2}^2\}}\right)$ such that
$$\mu F(x,s,t) \leq sF_s(x,s,t)+tF_t(x,s,t)+\sigma (s^2+t^2)\;\text{for all}\;x \in V,\;(s,t) \in \R^2,$$
where $\gamma_{2,i} = \mu_{\min}(\inf\limits_V V_i(x))^{\frac{1}{2}} $;\\
  Then problem (\ref{101a}) has infinitely many solutions $\{(u_k,v_k)\}_{k=1}^\infty$  such that (\ref{102a}) holds.}
\vskip2mm
\noindent
{\bf Remark 4.1.} From the proofs of Theorem 1.1 and Theorem 1.2, it is not difficult to see that Theorem 1.1 and Theorem 1.2 still hold even if one of $b_i, i=1,2$
is equal to zero, which implies that system (\ref{101}) is a coupling system between biharmonic-kirchhoff equation and biharmonic equation, and by comparing $(F_2), (F_4)$, $(F_5)$ with $(F_2)', (F_4)'$ and $(F_5)$, the biharmonic-kirchhoff equation plays a leading role rather than the biharmonic equation.

\vskip2mm
\noindent
{\bf Remark 4.2.} One can easily establish the  results corresponding to Theorem 4.1 and Theorem 4.2 for the following scalar biharmonic equation:
 \begin{eqnarray*}
 \Delta^2 u- a\Delta u + V(x)u=f(x,u),\;\;\;\;\hfill x\in V,
\end{eqnarray*}
by letting $u=v$, $a=a_1=a_2$, $V=V_1=V_2$ and $f(x,u)=F_u(x,u,v)=F_v(x,u,v)$ for all $x\in V$. At the moment,
$(F_1)$ and $(C_1)$ reduce to the following conditions, respectively:\\
$(f_1)$ there exists a function $a\in C(\R^+,\R^+)$ and a function $b:V\rightarrow \R^+$ with $b \in L^1(V)$ such that
$$|F(x,s)|, |f(x,s)| \leq a(|s|)b(x)$$
for all $x\in V$ and $s\in \R$, where $F(x,s)=\int_0^s f(x,t)dt$;\\
$(c_1)$  $\lim\limits_{|s|\rightarrow 0} \frac{a(|s|)}{|s|^2}=h<+\infty$.\\
In \cite{zhang2014}, in the Euclidean setting, the nonlinear term $f$ in system (\ref{104}) is required to satisfy that there exists constants $c_1,c_2>0$ and $p\in (2, 2_*)$ such that
 \begin{eqnarray}
 \label{b1}
|f(x,s)|\le c_1|s|+c_2|s|^{p-1}, \ \ \forall (x,s)\in \R^N\times \R
\end{eqnarray}
where $2_*=+\infty$ if $N\le 4$ and $2_*=\frac{2N}{N-4}$ if $N>4$. It is easy to verify that condition (\ref{b1}) excludes those examples with exponential growth such as
\begin{eqnarray} \label{b2}
f(x,s)=
   \begin{cases}
 a(x)e^s(s^4+4s^3),\;\;\;\;\hfill s>0,\\
 0,\;\;\;\;\hfill s=0,\\
  a(x)e^s(4s^3-s^4),\;\;\;\;\hfill s<0,\\
   \end{cases}
\end{eqnarray}
where $a(x) \in C(V,\R)$ and $0< \inf_V a(x) \leq \sup_V a(x)< 1 $. However, our assumptions $(f_1)$ and $(c_1)$ allow $f$ to be exponential growth and (\ref{b2}) satisfies all our assumptions in Theorem 4.2. The basic reason caused such difference is that $L^1(V)\subset L^\infty(V)$ and (\ref{221a}) holds in the locally finite graph setting, which in general dose not hold in the Euclidean setting $\R^N(N\ge 2)$.
\vskip3mm
 \noindent
\noindent{\bf Funding information}

\noindent
This work is supported by Yunnan Fundamental Research Projects (grant No: 202301AT070465).

\bibliographystyle{amsplain}

 \end{document}